\documentclass[review]{elsarticle}
\usepackage{amsmath}
\usepackage{amsthm}
\usepackage{caption}
\usepackage{lineno,hyperref}
\usepackage{subfigure}
\modulolinenumbers[5]

\newtheorem{lem}{Lemma}

\newtheorem{thm}{Theorem}
\newtheorem{cor}{Corollary}

\begin{document}

\begin{frontmatter}

\title{The crossing number of the complete 4-partite graph $K_{1,1,m,n}$}

%
%
%
\author[label1]{Xiwu Yang\corref{cor1}}
\ead{cinema@lnnu.edu.cn}
\author[label1]{Ni Lu}
\author[label1]{Xiaodong Chen}
\author[label2]{Yuansheng Yang}

\address[label1]{School of Mathematics, Liaoning Normal University, Dalian 116029,  P.R. China}
\address[label2]{Department of Computer Science, Dalian University of Technology, Dalian 116024, P.R.China}
\cortext[cor1]{Corresponding author}

\begin{abstract} Let $\textrm{cr}(G)$ denote the crossing number of a graph $G$.
The well-known Zarankiewicz's conjecture (ZC) asserted $\textrm{cr}(K_{m,n})$ in 1954.
In 1971, Harborth gave a conjecture (HC) on $\textrm{cr}(K_{x_1,...,x_n})$.
HC on $K_{1,m,n}$ is verified if ZC is true by Ho et al. in 2021.
In this paper, we showed the following results:
If both $m$ and $n$ are even, then
\[\textrm{cr}(K_{1,1,m,n})\geq \frac{1}{2}(\textrm{cr}(K_{m+1,n+3})+\textrm{cr}(K_{m+3,n+1})-mn-\frac{1}{4}(m^2+n^2));\]
If both $m$ and $n$ are odd, then
\[\textrm{cr}(K_{1,1,m,n})\geq \frac{1}{2}(\textrm{cr}(K_{1,m+1,n+1})+\textrm{cr}(K_{2,m,n})-\frac{1}{4}(m+1)(n+1)+1);\]
If $m$ is even and $n$ is odd, then
\begin{equation}\nonumber
\begin{split}
\textrm{cr}(K_{1,1,m,n})&\geq \frac{1}{4}(\textrm{cr}(K_{m+1,n+2})+\textrm{cr}(K_{m+3,n+2})+2\textrm{cr}(K_{2,m,n})
\\&-m(n+1)-\frac{1}{4}(n+1)^2).
\end{split}
\end{equation}
The lower bounds in our result imply that if both $m$ and $n$ are even and ZC is true, then HC on $K_{1,1,m,n}$ holds;
if at least one of $m$ and $n$ is odd and both ZC and HC on $K_{2,m,n}$ are true, then HC on $K_{1,1,m,n}$ holds.
\end{abstract}

\begin{keyword} Crossing number \sep Complete bipartite graph \sep Complete tripartite graph \sep Harborth's conjecture \sep Zarankiewicz's conjecture
\MSC[2010] 05C10
\end{keyword}

\end{frontmatter}


\section{Introduction}

A drawing of a graph $G$ is said to be a \emph{good} drawing, provided that no edges crosses itself,
no adjacent edges cross each other, no two edges cross more than once,
and no three edges cross at a common point.
In a good drawing, a common point of two edges other than endpoints is a \emph{crossing}.
The crossing number of $G$, denoted by $\textrm{cr}(G)$, is the minimum number of
crossings among all good drawings of $G$ in the plane.

Zarankiewicz\cite{Zar} introduced a good drawing of the complete bipartite graph $K_{m,n}$,
in which $\textrm{cr}(K_{m,n})\leq Z(m,n)$,
where $Z(m,n)=\lfloor\frac{m}{2}\rfloor\lfloor\frac{m-1}{2}\rfloor\lfloor\frac{n}{2}\rfloor\lfloor\frac{n-1}{2}\rfloor$
and $\lfloor x \rfloor$ denotes the largest integer no more than $x$.
Zarankiewicz \cite{Zar} also claimed $\textrm{cr}(K_{m,n})=Z(m,n)$,
but there is a subtle error in his proof, found by Kainen and Ringel respectively,
as described in \cite{RefGuy}.
Thus, the following equality is known as Zarankiewicz's conjecture (ZC):
\begin{equation}\label{A0}
\textrm{cr}(K_{m,n})=Z(m,n).
\end{equation}
Kleitman \cite{Kle} verified that \eqref{A0} is true for $\textrm{min}(m,n)\leq 6$.
Christian et al. \cite{Chr} proved that
for every $m$ if there is an integer $N(m)$ such that Zarankiewicz's conjecture is true for $K_{m,n}$ with $n \leq N(m)$,
then it is true for each positive integer $n$.

Harborth \cite{Har} showed an upper bound on the crossing number of the complete $n$-partite graph $K_{x_1,...,x_n}$ as follows.
\begin{equation}\label{A1}
\begin{split}
\textrm{cr}(K_{x_1,...,x_n})\leq&\frac{1}{8}\left(\sum\limits_{1\leq i<j<r<s\leq n}3x_ix_jx_rx_s+3\tbinom{\lfloor c/2\rfloor}{2}\right.
\\&\phantom{=\;\;}\left.-\sum\limits_{1\leq i<j\leq n}\lfloor\frac{c-((x_i\ \mbox{mod}\ 2)+(x_j\ \mbox{mod}\ 2))}{2}\rfloor x_ix_j\right)
\\&+\sum\limits_{i=1}^n \lfloor\frac{x_i}{2}\rfloor\lfloor\frac{x_i-1}{2}\rfloor\lfloor\frac{m-x_i}{2}\rfloor\lfloor\frac{m-x_i-1}{2}\rfloor
\\&-\sum\limits_{1\leq i\leq j\leq n} \lfloor\frac{x_i}{2}\rfloor\lfloor\frac{x_i-1}{2}\rfloor\lfloor\frac{x_j}{2}\rfloor\lfloor\frac{x_j-1}{2}\rfloor,
\end{split}
\end{equation}
where $m=\sum_{i=1}^nx_i$ and $c$ is the number of odd number $x_i$ for $1\leq i\leq n$.
Harborth \cite{Har} also conjectured that equality holds in \eqref{A1}.
We call it Harborth's conjecture (HC) in this paper.
For $n=2$, HC on $K_{m,n}$ is equivalent to \eqref{A0}.
HC on $K_{1,m,n}$ is equivalent to
\begin{equation}\label{A2}
\textrm{cr}(K_{1,m,n})=Z(m+1,n+1)-\lfloor\frac{m}{2}\rfloor\lfloor\frac{n}{2}\rfloor.
\end{equation}
HC on $K_{2,m,n}$ is equivalent to
\begin{equation}\label{A3}
\textrm{cr}(K_{2,m,n})=Z(m+2,n+2)-mn.
\end{equation}
HC on $K_{1,1,m,n}$ is equivalent to
\begin{equation}\label{A5}
\textrm{cr}(K_{1,1,m,n})=Z(m+2,n+2)-mn+\lfloor\frac{m}{2}\rfloor\lfloor\frac{n}{2}\rfloor.
\end{equation}

Asano \cite{Asa} proved that both \eqref{A2} and \eqref{A3} are true for $m=3$.
Huang and Zhao \cite{Hua} showed that \eqref{A2} is true for $m=4$.
Ho \cite{Ho08} obtained that
\begin{equation}\label{A10}
\textrm{cr}(K_{1,m,n})\geq \textrm{cr}(K_{m+1,n+1})-\left\lfloor\frac{n}{m}\lfloor\frac{m}{2}\rfloor\lfloor\frac{m+1}{2}\rfloor\right\rfloor;
\end{equation}
and if $m$ is even, then
\begin{equation}\label{A4}
\textrm{cr}(K_{1,m,n})\geq \frac{1}{2}(\textrm{cr}(K_{m+1,n+2})+\textrm{cr}(K_{m+1,n})-\frac{m}{2}(\frac{m}{2}+n-1)).
\end{equation}
The first author and Wang \cite{Yan} asserted that if both $m$ and $n$ are odd, then
\begin{equation}\label{A6}
\textrm{cr}(K_{1,m,n})\geq \frac{1}{2}(\textrm{cr}(K_{m+2,n})+\textrm{cr}(K_{m,n+2})-\lfloor\frac{m}{2}\rfloor^2-\lfloor\frac{n}{2}\rfloor^2).
\end{equation}
Suppose ZC is true. Then \eqref{A10} implies that \eqref{A2} is true when both $m$ and $n$ are even;
\eqref{A4} implies that \eqref{A2} holds when $m$ is even and $n$ is odd; \eqref{A6} implies that \eqref{A2} is true if both $m$ and $n$ are odd.
Ho \cite{Ho11,Ho13} showed that \eqref{A3} is true for $m=4$ and \eqref{A5} is true for $m=3$.
Ho \cite{Ho09} obtained the crossing numbers of $K_{1,1,1,1,n}$, $K_{1,1,1,2,n}$ and $K_{1,2,2,n}$,
in which Ho verified that HC on each of $K_{1,1,1,1,n}$, $K_{1,1,1,2,n}$ and $K_{1,2,2,n}$ is true.
For other results about the crossing number of graphs, readers can refer to \cite{Cla,Sch18, Sch17}.

In this paper, we obtained the lower bounds on $\textrm{cr}(K_{1,1,m,n})$ on all different parities of $m$ and $n$.
As corollaries, we gave the exactly values of $\textrm{cr}(K_{1,1,4,4})$ and $\textrm{cr}(K_{1,1,3,n})$,
and proved that HC on $K_{1,1,m,n}$ holds if both ZC and HC on $K_{2,m,n}$ are true.

\section{Terminology and Drawing Lemma}

For some terminology and notation not defined here, readers can refer to \cite{Bon}.
In this paper, we only consider simple graphs.
Let $G$ be a graph with $A,B\subseteq E(G)$.
In a good drawing $D$ of $G$, the number of crossings crossed by one edge in $A$
and the other edge in $B$ is denoted by $\textrm{cr}_D(A,B)$.
Especially, $\textrm{cr}_D(A,A)$ is denoted by $\textrm{cr}_D(A)$ for short.
The number of crossings of $D$, denoted by $\textrm{cr}(D)$, is $\textrm{cr}_D(E(G))$.
The following two equalities are trivial.
\begin{eqnarray}
\textrm{cr}_D(A\cup B)&=&\textrm{cr}_D(A)+\textrm{cr}_D(B)+\textrm{cr}_D(A,B),\label{I1}\\
\textrm{cr}_D(A,B\cup C)&=&\textrm{cr}_D(A,B)+\textrm{cr}_D(A,C),\label{I2}
\end{eqnarray}
where $A$, $B$ and $C$ are pairwise disjoint subsets of $E(G)$. %

For a vertex $v$ in $G$, let $\textrm{cr}_D(v)$ denote the number of crossings in all the edges incident with $v$.
For a vertex $u$ in $G$ distinct from $v$, let $\textrm{cr}_D(u,v)$ denote the number of crossings crossed by one edge incident with $u$ and the other edge incident with $v$.
The \emph{neighborhood} of $v$, denoted by $N(v)$, is the set of vertices adjacent to $v$.
The set of edges incident to $v$ is denoted by $E(v)$.
For a given $D$, $E(v)$ induces a natural cyclic order $\pi$ of $N(v)$, which is called \emph{rotation}, denoted by $\pi_D(v)$ of $v$ in $D$.
For $U\subset N(v)$, let
$\pi_D^U(v)$ denote the subrotation of $\pi_D(v)$ from $U$.
For $u,w\in N(v)$, let $\pi_D(vu,vw)$ denote the sublist in $\pi_D(v)$ from the vertice following $u$ to the vertice before $w$.
For $V_1,V_2\subseteq V(G)$, let $E(V_1,V_2)=\{v_1v_2|v_1\in V_1, v_2\in V_2\}$.
For $x\notin V(G)$ and $v\in V(G)$, let $G^{xTv}$ be the graph such that $V(G^{xTv}) =V(G)\cup \{x\}$ and
$E(G^{xTv}) =E(G)\cup \{xu|u\in N(v)\}$.

The following lemma is generalization of the good drawings in \cite{Ho08} and \cite{Hua}.

\begin{lem}\label{lem1}
Let $G$ be a connected graph, and $v$ be a vertex in $G$ with degree $p+q$ and $N(v)=U\cup W$, where $U=\{u_0,...,u_{p-1}\}$ and $W=\{w_0,...,w_{q-1}\}$.
Assume $G_2=G_1^{yTx}$, where $V(G_1)=V(G)\cup \{x\}$ and $E(G_1)=E(G)\cup \{xu_i|0\leq i\leq p-1\}\cup \{xv\}-\{vu_j|0\leq j\leq p-1\}$.
Suppose $D$ is a good drawing of $G$ with
$\pi_D^U(v)=(u_0,...,u_{p-1})$ and
$W_i=\pi_D(vu_i,vu_{i+1})$ for $0\leq i\leq p-1$.
Then $q=\sum_{i=0}^{p-1}|W_i|$.
If $p$ is even, then we can construct a good drawing $D^1_k$ of $G_1$ with $\textrm{cr}(D^1_k)=\textrm{cr}(D)+\sum_{s=k}^{k+\frac{p}{2}-2}(k+\frac{p}{2}-s-1)|W_s|+\sum_{s=k+\frac{p}{2}}^{k+p-1}(s+1-k-\frac{p}{2})|W_s|$,
and a good drawing $D^2_k$ of $G_2$ with $\textrm{cr}(D^2_k)=\textrm{cr}(D)+\textrm{cr}_D(E(\{v\},U),E(G)-E(v))+\frac{p}{2}(q+\frac{p}{2}-1)$, for $0\leq k\leq p-1$.
\end{lem}

\begin{proof} By the definition of a good drawing, let $N(v,\varepsilon)=\{s\in R^2:||s-v||<\varepsilon\}$, where $R^2$ denotes the plane and $\varepsilon$ is a sufficiently small positive number
such that both $V(G)-\{v\}$ and $E(G)-E(v)$ are located outside $N(v,\varepsilon)$.
See Fig.\ref{F1}(a) for $|U|=6$, for example.

For $0\leq k \leq p-1$,
we can construct a good drawing $D^1_k$ of $G_1$ by modifying $D$ in the following steps:

Step $(i)$. Draw a vertex $x$ on a point of $vu_k$ in $N(v,\varepsilon)$ such that $vu_k$ is split into edges $vx$ and $xu_k$.

Step $(ii)$. For $k+1\leq i\leq k+\frac{p}{2}-1$,
draw an edge $u_ix$ next to the edge $u_iv$ such that $u_ix$ crosses each edge which $u_iv$ crosses;
Moreover, let $u_ix$ crosses each $vw_j$ and each $vu_r$ for $w_j\in \cup_{s=k}^{i-1} W_s$ (mod $p$ for the subscript $s$ of $W_s$) and $k+1\leq r\leq i-1$ respectively.
For $k+\frac{p}{2}\leq i\leq k+p-1$,
draw an edge $u_ix$ next to the edge $u_iv$ such that $u_ix$ crosses each edge which $u_iv$ crosses;
Moreover, let $u_ix$ crosses each $vw_j$ and each $vu_r$ for $w_j\in \cup_{s=i}^{k+p-1} W_s$ and $i+1\leq r\leq k+p-1$ respectively.

Step $(iii)$. Remove each edge $vu_i$ for $0\leq i\leq p-1$ and $i\neq k$.

See Fig.\ref{F1}(b) for $|U|=6$, for example.
By the definition of a good drawing, $\textrm{cr}_D(E(\{v\},U))=0$.
By \eqref{I1} and \eqref{I2}, we have
\begin{equation}\label{C1}
\textrm{cr}(D)=\textrm{cr}_D(E(G)-E(\{v\},U))+\textrm{cr}_D(E(G)-E(\{v\},U),E(\{v\},U)).
\end{equation}
For $0\leq k \leq p-1$, by the same argument, we assert $\textrm{cr}_{D^1_k}(E(x))=0$ and
\begin{equation}\label{C2}
\textrm{cr}(D^1_k)=\textrm{cr}_{D^1_k}(E(G_1)-E(x))+\textrm{cr}_{D^1_k}(x).
\end{equation}
For $0\leq k \leq p-1$, it is easy to obtain that the number of $\textrm{cr}_{D^1_k}(x)$ in $N(v,\varepsilon)$ is $\sum_{s=k}^{k+\frac{p}{2}-2}(k+\frac{p}{2}-s-1)|W_s|+\sum_{s=k+\frac{p}{2}}^{k+p-1}(s+1-k-\frac{p}{2})|W_s|$,
and the number of $\textrm{cr}_{D^1_k}(x)$ outside $N(v,\varepsilon)$ is $\textrm{cr}_D(E(\{v\},U),E(G)-E(\{v\},U))$.
Hence, for $0\leq k \leq p-1$ we claim
\begin{equation}\label{C3}
\begin{split}
\textrm{cr}_{D^1_k}(x)=&\sum_{s=k}^{k+\frac{p}{2}-2}(k+\frac{p}{2}-s-1)|W_s|+\sum_{s=k+\frac{p}{2}}^{k+p-1}(s+1-k-\frac{p}{2})|W_s|\\&+\textrm{cr}_D(E(G)-E(\{v\},U),E(\{v\},U)).
\end{split}
\end{equation}
For $0\leq k \leq p-1$, since $\textrm{cr}_{D^1_k}(E(G_1)-E(x))=\textrm{cr}_D(E(G)-E(\{v\},U))$, \eqref{C1}-\eqref{C3} imply
\begin{equation}\label{Z1}
\textrm{cr}(D^1_k)=\textrm{cr}(D)+\sum_{s=k}^{k+\frac{p}{2}-2}(k+\frac{p}{2}-s-1)|W_s|+\sum_{s=k+\frac{p}{2}}^{k+p-1}(s+1-k-\frac{p}{2})|W_s|.
\end{equation}

For $0\leq k \leq p-1$,
we can construct a good drawing $D^2_k$ of $G_2$ by modifying $D$ in the following steps:

Step $(i)$. Do the first two steps in drawing $D^1_{k+\frac{p}{2}}$.

Step $(ii)$. Draw a vertex $y$ on a point of $vu_k$ in $N(v,\varepsilon)$ such that $vu_k$ is split into edges $vy$ and $yu_k$.
Draw an edge $vu_{k+\frac{p}{2}}$ next to the edge $xu_{k+\frac{p}{2}}$ such that
$vu_{k+\frac{p}{2}}$ crosses each $xu_r$ for $k+\frac{p}{2}+1\leq r\leq k+p-1$ and
each edge which $xu_{k+\frac{p}{2}}$ crosses.

Step $(iii)$. For $k+1\leq i\leq k+\frac{p}{2}-1$,
draw an edge $u_iy$ next to the edge $u_iv$ such that $u_iy$ crosses each edge which $u_iv$ crosses;
Moreover, let $u_iy$ crosses each of $vw_j$,$vu_r$ and $xu_t$ for $w_j\in \cup_{s=k}^{i-1} W_s$, $k+1\leq r\leq i-1$ and $k\leq t\leq i-1$ respectively.
For $k+\frac{p}{2}\leq i\leq k+p-1$,
draw an edge $u_iy$ next to the edge $u_iv$ such that $u_iy$ crosses each edge which $u_iv$ crosses;
Moreover, let $u_iy$ crosses each of $vw_j$, $vu_r$ and $xu_r$ for $w_j\in \cup_{s=i}^{k+p-1} W_s$ and $i+1\leq r\leq k+p-1$ respectively.

Step $(iv)$. Remove each edge $vu_i$ for $0\leq i\leq p-1$ and $i\neq k$.

See Fig.\ref{F1}(c) for $|U|=6$, for example.
For $0\leq k \leq p-1$, it is easy to check that $\textrm{cr}(D^2_k)=\textrm{cr}(D^1_{k+\frac{p}{2}})+\textrm{cr}_{D^2_k}(y)$.
For $0\leq k \leq p-1$, the number of $\textrm{cr}_{D^2_k}(y)$ in $N(v,\varepsilon)$ is the sum of $\textrm{cr}_{D^2_k}(x,y)$ and $\textrm{cr}_{D^2_k}(v,y)$,
which are $\frac{p}{2}(\frac{p}{2}-1)$ and $\sum_{s=k}^{k+\frac{p}{2}-2}(k+\frac{p}{2}-s-1)|W_s|+\sum_{s=k+\frac{p}{2}}^{k+p-1}(s+1-k-\frac{p}{2})|W_s|$ respectively.
For $0\leq k \leq p-1$, the number of $\textrm{cr}_{D^2_k}(y)$ outside $N(v,\varepsilon)$ is $\textrm{cr}_D(E(\{v\},U),E(G)-E(v))$.
Hence,
\begin{equation}\label{Z0}
\begin{split}
\textrm{cr}_{D^2_k}(y)=&\sum_{s=k}^{k+\frac{p}{2}-2}(k+\frac{p}{2}-s-1)|W_s|+\sum_{s=k+\frac{p}{2}}^{k+p-2}(s+1-k-\frac{p}{2})|W_s|\\&+\frac{p}{2}|W_{k-1}|+\frac{p}{2}(\frac{p}{2}-1)+\textrm{cr}_D(E(\{v\},U),E(G)-E(v)).
\end{split}
\end{equation}
\captionsetup[figure]{labelfont={bf},name={Fig.},labelsep=period}
\begin{figure}
  \centering
  \subfigure[$N(v,\varepsilon)$ in $D$.]{
    \label{fig:subfig:a} 
    \includegraphics[width=0.3\textwidth]{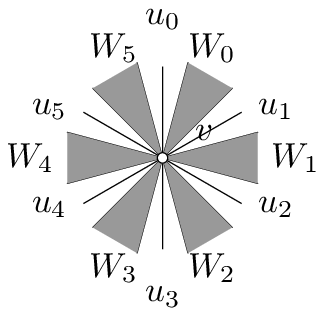}}
  \subfigure[$N(v,\varepsilon)$ in $D^1_3$.]{
    \label{fig:subfig:a} 
    \includegraphics[width=0.3\textwidth]{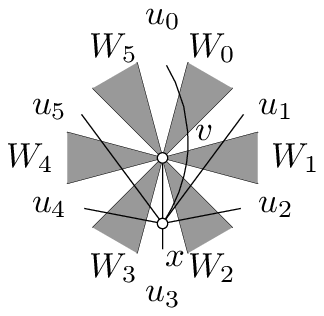}}
  \subfigure[$N(v,\varepsilon)$ in $D^2_0$.]{
    \label{fig:subfig:b} 
    \includegraphics[width=0.3\textwidth]{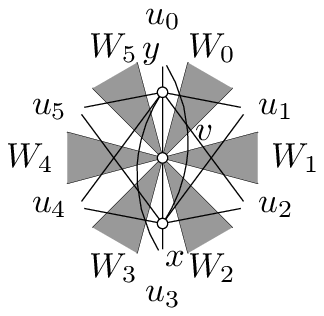}}
  \caption{Subdrawings of $D$, $D^1_3$ and $D^2_0$ for $|U|=6$. Each dark triangle denotes the part of $E(\{v\},W_i)$ in $N(v,\varepsilon)$ for $0\leq i\leq 5$.
  Each semi-edge represents the part of $vu_i$, $xu_i$ or $yu_i$ in $N(v,\varepsilon)$ for $0\leq i\leq 5$. }\label{F1}
\end{figure}
By \eqref{Z1}, we have
\begin{equation}\label{Z2}
\textrm{cr}(D^1_{k+\frac{p}{2}})=\textrm{cr}(D)+\sum_{s=k+\frac{p}{2}}^{k+p-2}(k+p-s-1)|W_s|+\sum_{s=k+p}^{k+p+\frac{p}{2}-1}(s+1-k-p)|W_s|.
\end{equation}
Since $\sum_{s=k+p}^{k+p+\frac{p}{2}-1}(s+1-k-p)|W_s|=\sum_{s=k}^{k+\frac{p}{2}-1}(s+1-k)|W_s|$,
\eqref{Z2} implies
\begin{equation}\label{Z3}
\begin{split}
\textrm{cr}(D^1_{k+\frac{p}{2}})=&\sum_{s=k+\frac{p}{2}}^{k+p-2}(k+p-s-1)|W_s|+\sum_{s=k}^{k+\frac{p}{2}-2}(s+1-k)|W_s|\\&+\frac{p}{2}|W_{k+\frac{p}{2}-1}|+\textrm{cr}(D).
\end{split}
\end{equation}
Recall that $\textrm{cr}(D^2_k)=\textrm{cr}(D^1_{k+\frac{p}{2}})+\textrm{cr}_{D^2_k}(y)$.
Moreover, by \eqref{Z0}, \eqref{Z3} and $q=\sum_{i=0}^{p-1}|W_i|$, we have $\textrm{cr}(D^2_k)=\textrm{cr}(D)+\textrm{cr}_D(E(\{v\},U),E(G)-E(v))+\frac{p}{2}(q+\frac{p}{2}-1)$ for $0\leq k \leq p-1$.
\end{proof}

\section{Lower bounds on $\textrm{cr}(K_{1,1,m,n})$}

In the following proof, let $O$, $X$, $Y$ and $Z$ be the independent sets of $K_{1,1,m,n}$,
where $O=\{o\}$, $X=\{x\}$, $Y=\{y_1,\ldots,y_m\}$ and $Z=\{z_1,\ldots,z_n\}$.
By the definition of a good drawing, \eqref{I1} and \eqref{I2}, we have

\begin{lem}\label{lem2}
Let $D$ be an arbitrary good drawing of $K_{1,1,m,n}$. Then
$\textrm{cr}(D)=\textrm{cr}_{D}(E(Y,Z))+\textrm{cr}_{D}(E(O,X), E(Y,Z))+\textrm{cr}_{D}(E(X,Y),E(Y,Z))+\textrm{cr}_{D}(E(O,Z),$
$E(X,Y)\cup E(Y,Z))+\textrm{cr}_{D}(E(O,Y), E(X,Z))+\textrm{cr}_{D}(E(O,Y), E(X,Y)\cup E(Y,Z))
+\textrm{cr}_{D}(E(X,Z), E(O,Z)\cup E(Y,Z)).
$
\end{lem}

\begin{lem}\label{lem3}
Let $D$ be a good drawing of $K_{1,1,m,n}$. Then
$\textrm{cr}_{D}(E(O,X),E(Y,Z))=\textrm{cr}_{D}(E(O,X),E(K_{1,1,m,n})-E(O,X))\leq \textrm{cr}(D)-\textrm{cr}(K_{2,m,n}).$
\end{lem}

\begin{proof}By Lemma \ref{lem2} and \eqref{I1}, we assert
$\textrm{cr}(D)=\textrm{cr}_{D}(E(O,X))+\textrm{cr}_{D}(E(K_{1,1,m,n})-E(O,X))+\textrm{cr}_{D}(E(O,X),E(K_{1,1,m,n})-E(O,X)).$
By the definition of a good drawing, $\textrm{cr}_{D}(E(O,X))=0$.
Since the graph induced by $E(K_{1,1,m,n})-E(O,X)$ is isomorphic to $K_{2,m,n}$, we claim $\textrm{cr}_{D}(E(K_{1,1,m,n})-E(O,X))\geq \textrm{cr}(K_{2,m,n})$.
Hence, we have $\textrm{cr}_{D}(E(O,X),E(K_{1,1,m,n})-E(O,X))\leq \textrm{cr}(D)-\textrm{cr}(K_{2,m,n})$.
Since $\textrm{cr}_{D}(E(O,X),E(O,Y)\cup E(O,Z)\cup E(X,Y)\cup E(X,Z)=0$ by the definition of a good drawing, 
we have $\textrm{cr}_{D}(E(O,X),E(K_{1,1,m,n})-E(O,X))=\textrm{cr}_{D}(E(O,X),E(Y,Z))$.
\end{proof}

\begin{thm}\label{thm1}
If both $m$ and $n$ are even, then
\[\textrm{cr}(K_{1,1,m,n})\geq \frac{1}{2}(\textrm{cr}(K_{m+1,n+3})+\textrm{cr}(K_{m+3,n+1})-mn-\frac{1}{4}(m^2+n^2)).\]
\end{thm}

\begin{proof}
Let $D$ be a good drawing of $K_{1,1,m,n}$.
Assume $G_1$ is the graph obtained by deleting $E(X,Z)$ from $K_{1,1,m,n}$, and $D_1$ is the drawing of $G_1$ obtained by deleting $E(X,Z)$ from $D$.
By Lemma \ref{lem2}, we claim
\begin{equation}\label{B4}
\begin{split}
\textrm{cr}(D_1)=\textrm{cr}(D)-\textrm{cr}_{D}(E(X,Z), E(O,Y)\cup E(O,Z)\cup E(Y,Z)).
\end{split}
\end{equation}
Suppose $G_3=G_2^{z_{n+2}Tz_{n+1}}$,
where $V(G_2)=V(G_1)\cup \{z_{n+1}\}$ and $E(G_2)=E(G_1)\cup \{z_{n+1}y_i|1\leq i\leq m\}\cup \{z_{n+1}o\}-\{oy_i|1\leq i\leq m\}$.
It is easy to check that $G_3$ is isomorphic to $K_{m+1,n+3}$ with the independent sets $O\cup Y$ and $X\cup Z\cup\{z_{n+1},z_{n+2}\}$.

In $D_1$, we give $x$ another label $z_0$.
By renaming the vertices of $y_i$ if necessary,
we assume $\pi_{D_1}^Y(o)=(y_1,...,y_m)$ and
$Z_i=\pi_D(oy_i,oy_{i+1})$ for $1\leq i\leq m$.
Then $n+1=\sum_{i=1}^m|Z_i|$.
Since $m$ is even, we can construct a good drawing $D_2$ of $G_3$ with $\textrm{cr}(D_1)+\textrm{cr}_{D_1}(E(O,Y),E(G_1)-E(o))+\frac{m}{2}(n+\frac{m}{2})$ crossings by Lemma \ref{lem1}.
For example, a good drawing $D$ of $K_{1,1,4,4}$ is illustrated in Fig.\ref{F2}.
A good drawing $D_2$ of $K_{5,7}$ constructed from $D$ in Fig.\ref{F2} is illustrated in Fig.\ref{F3}.
Hence, $\textrm{cr}(D_2)\geq \textrm{cr}(K_{m+1,n+3})$.
Since $\textrm{cr}_{D_1}(E(O,Y),E(G_1)-E(o))=\textrm{cr}_{D}(E(O,Y),E(X,Y)\cup E(Y,Z))$, we have
\begin{equation}\label{B5}
\begin{split}
\textrm{cr}(D_2)=\textrm{cr}(D_1)+\textrm{cr}_{D}(E(O,Y),E(X,Y)\cup E(Y,Z))+\frac{m}{2}(n+\frac{m}{2}).
\end{split}
\end{equation}
\captionsetup[figure]{labelfont={bf},name={Fig.},labelsep=period}
\begin{figure}[htp]
	\centering
	\begin{minipage}[b]{1.0\textwidth}
		\centering
		\includegraphics[width= 1.0\textwidth]{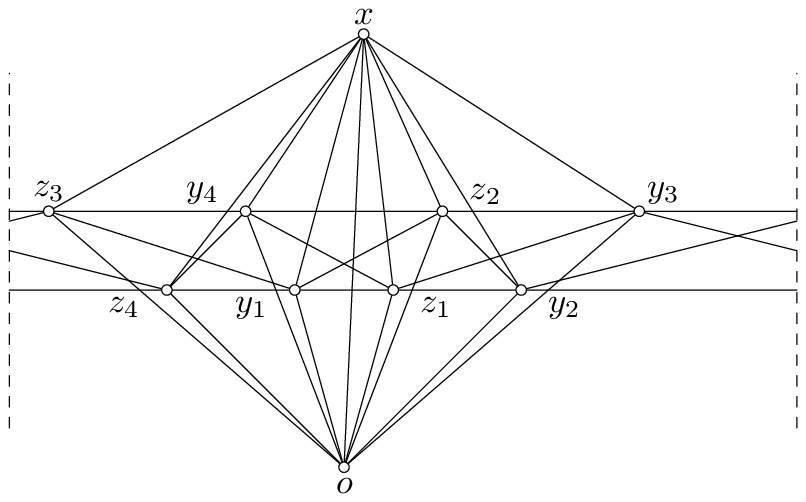}
		\begin{minipage}[t]{\textwidth} \caption{A good drawing $D$ of $K_{1,1,4,4}$ in a cylinder.}\label{F2} \end{minipage}
	\end{minipage}
	\hfill
	\begin{minipage}[b]{1.0\textwidth}
		\centering
		\includegraphics[width= 1.0\textwidth]{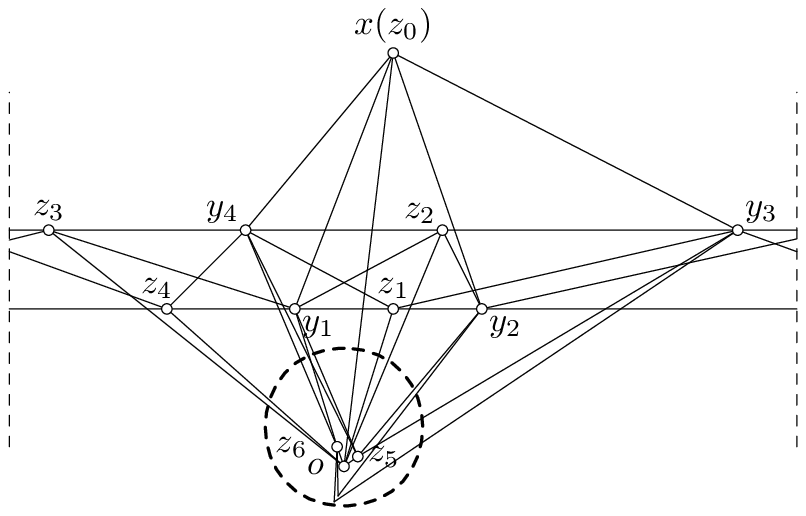}
		\begin{minipage}[t]{\textwidth} \caption{A good drawing $D_2$ of $K_{5,7}$ constructed from $D$ in Fig.\ref{F2}. Several locations of vertices were adjusted for visual clarity.
  The dashed circle denotes the boundary of $N(o,\varepsilon)$ in $D$, which is not a part of $D_2$.}\label{F3} \end{minipage}
	\end{minipage}
\end{figure}
Similarly, let $H_1=K_{1,1,m,n}-E(O,Y)$ and $D_3$ be the drawing of $H_1$ obtained by deleting $E(O,Y)$ from $D$.
By Lemma \ref{lem2}, we assert
\begin{equation}\label{B2}
\begin{split}
\textrm{cr}(D_3)=\textrm{cr}(D)-\textrm{cr}_{D}(E(O,Y), E(X,Y)\cup E(X,Z)\cup E(Y,Z)).
\end{split}
\end{equation}
Assume $H_3=H_2^{y_{m+2}Ty_{m+1}}$,
where $V(H_2)=V(H_1)\cup \{y_{m+1}\}$ and $E(H_2)=E(H_1)\cup \{y_{m+1}z_i|1\leq i\leq n\}\cup \{y_{m+1}x\}-\{xz_i|1\leq i\leq n\}$.
It is easy to check that $H_3$ is isomorphic to $K_{m+3,n+1}$ with the independent sets $O\cup Y\cup\{y_{m+1},y_{m+2}\}$ and $X\cup Z$.

In $D_3$, we give $o$ another label $y_0$.
By renaming the vertices of $z_i$ if necessary,
we assume $\pi_{D_3}^Z(x)=(z_1,...,z_n)$ and
$Y_i=\pi_D(oz_i,oz_{i+1})$ for $1\leq i\leq n$.
Then $m+1=\sum_{i=1}^n|Y_i|$.
Since $n$ is even, we can construct a good drawing $D_4$ of $H_3$ with $\textrm{cr}(D_3)+\textrm{cr}_{D_3}(E(X,Z),E(H_1)-E(x))+\frac{n}{2}(m+\frac{n}{2})$ crossings by Lemma \ref{lem1}.
For example, a good drawing $D_4$ of $K_{5,7}$ constructed from $D$ in Fig.\ref{F2} is illustrated in Fig.\ref{F4}.
Hence, $\textrm{cr}(D_4)\geq \textrm{cr}(K_{m+3,n+1})$.
Since $\textrm{cr}_{D_3}(E(X,Z),E(H_1)-E(x))=\textrm{cr}_{D}(E(X,Z),E(O,Z)\cup E(Y,Z))$, we claim
\begin{equation}\label{B3}
\begin{split}
\textrm{cr}(D_4)=\textrm{cr}(D_3)+\textrm{cr}_{D}(E(X,Z),E(O,Z)\cup E(Y,Z))+\frac{n}{2}(m+\frac{n}{2}).
\end{split}
\end{equation}

By \eqref{B4}-\eqref{B3}, we assert
\begin{equation}\label{B6}
\begin{split}
\textrm{cr}(D_2)+\textrm{cr}(D_4)=&2\textrm{cr}(D)-2\textrm{cr}_{D}(E(O,Y),E(X,Z))\\&+\frac{n}{2}(m+\frac{n}{2})+\frac{m}{2}(n+\frac{m}{2}).
\end{split}
\end{equation}
By the definition of a good drawing, $\textrm{cr}_{D}(E(O,Y),E(X,Z))\geq 0$.
\eqref{B6} implies $\textrm{cr}(D)\geq \frac{1}{2}(\textrm{cr}(K_{m+1,n+3})+\textrm{cr}(K_{m+3,n+1})-mn-\frac{1}{4}(m^2+n^2))$,
by $\textrm{cr}(D_2)\geq \textrm{cr}(K_{m+1,n+3})$, $\textrm{cr}(D_4)\geq \textrm{cr}(K_{m+3,n+1})$ and  $\textrm{cr}_{D}(E(O,Y),E(X,Z))\geq 0$.
By the definition of crossing number, we have $\textrm{cr}(K_{1,1,m,n})\geq\frac{1}{2}(\textrm{cr}(K_{m+1,n+3})+\textrm{cr}(K_{m+3,n+1})-mn-\frac{1}{4}(m^2+n^2))$.
\end{proof}
\captionsetup[figure]{labelfont={bf},name={Fig.},labelsep=period}
\begin{figure}[htp]
	\centering
	\begin{minipage}[b]{1.0\textwidth}
		\centering
		\includegraphics[width= 1.0\textwidth]{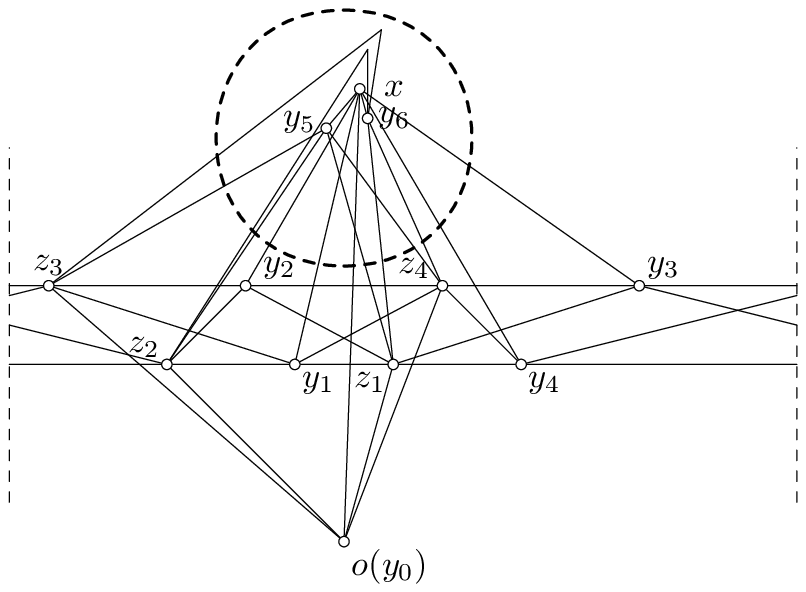}
		\begin{minipage}[t]{\textwidth} \caption{A good drawing $D_4$ of $K_{5,7}$ constructed from $D$ in Fig.\ref{F2}.
  Several locations of vertices were adjusted for visual clarity.
  The dashed circle denotes the boundary of $N(x,\varepsilon)$ in $D$, which is not a part of $D_4$.}\label{F4} \end{minipage}
	\end{minipage}
\end{figure}

\begin{cor}\label{cor1}
If both $m$ and $n$ are even and ZC is true, then HC on $K_{1,1,m,n}$ is true.
\end{cor}

\begin{proof} By \eqref{A1}, it suffices to prove $\textrm{cr}(K_{1,1,m,n})\geq Z(m+2,n+2)-mn+\lfloor\frac{m}{2}\rfloor\lfloor\frac{n}{2}\rfloor.$
Since both $m$ and $n$ are even, let $m=2r$ and $n=2s$.
By Theorem \ref{thm1}, $\textrm{cr}(K_{1,1,2r,2s})\geq \frac{1}{2}(\textrm{cr}(K_{2r+1,2s+3})+\textrm{cr}(K_{2r+3,2s+1})-4rs-\frac{1}{4}(4r^2+4s^2))$.
Since ZC is true, $\textrm{cr}(K_{2r+1,2s+3})=r^2(s+1)^2$ and $\textrm{cr}(K_{2r+3,2s+1})=(r+1)^2s^2$.
Hence, $\textrm{cr}(K_{1,1,2r,2s})\geq r^2s^2+sr(s+r-2)=Z(m+2,n+2)-mn+\lfloor\frac{m}{2}\rfloor\lfloor\frac{n}{2}\rfloor$.
\end{proof}

\begin{lem}[\cite{Kle}]\label{lem4}
$\textrm{cr}(K_{5,n})=Z(5,n)$.
\end{lem}

\begin{cor}\label{cor2}
 $\textrm{cr}(K_{1,1,4,4})=24.$
\end{cor}

\begin{proof} By Lemma \ref{lem4}, $\textrm{cr}(K_{5,7})=36$.
Therefore, $\textrm{cr}(K_{1,1,4,4})= 24$ by the same arguments in the proof of Corollary \ref{cor1}.
\end{proof}

\begin{thm}\label{thm2}
If both $m$ and $n$ are odd, then
\[\textrm{cr}(K_{1,1,m,n})\geq \frac{1}{2}(\textrm{cr}(K_{1,m+1,n+1})+\textrm{cr}(K_{2,m,n})-\frac{1}{4}(m+1)(n+1)+1).\]
\end{thm}

\begin{proof}Assume $V(G_1)=V(K_{1,1,m,n})\cup \{y_0\}$ and $E(G_1)=E(K_{1,1,m,n})\cup \{y_0z_i|$ $1\leq i\leq n\}\cup \{y_0o,y_0x\}-\{oz_i|1\leq i\leq n\}\cup \{ox\}$.
Let $D$ be a good drawing of $K_{1,1,m,n}$.
In $D$, we give $x$ another label $z_0$.
By renaming the vertices of $z_i$ if necessary,
we assume $\pi_D^{X\cup Z}(o)=(z_0,...,z_n)$ and
$Y_i=\pi_D(oz_i,oz_{i+1})$ for $0\leq i\leq n$.
Then $m=\sum_{i=0}^n|Y_i|$.
Let $N(o,\varepsilon)=\{s\in R^2:||s-o||<\varepsilon\}$, where $R^2$ denotes the plane and $\varepsilon$ is a sufficiently small positive number
such that both $X\cup Y\cup Z$ and $E(X,Y)\cup E(X,Z)\cup E(Y,Z)$ are located outside $N(o,\varepsilon)$.
Since $n+1$ is even, by the proof of Lemma \ref{lem1}, we can get two good drawings $D^1_0$ and $D^1_{\frac{n+1}{2}}$ of $G_1$ with
\begin{eqnarray}
\begin{split}
\textrm{cr}(D^1_0)=&\textrm{cr}(D)+\sum_{i=0}^{(n+1)/2-2}((n+1)/2-i-1)|Y_i|\\&+\sum_{i=(n+1)/2}^{n}(i+1-(n+1)/2)|Y_i|,
\end{split}\label{B7}\\
\textrm{cr}(D^1_{\frac{n+1}{2}})=\textrm{cr}(D)+\sum_{i=(n+1)/2}^{n-1}(n-i)|Y_i|+\sum_{i=0}^{(n+1)/2-1}(i+1)|Y_i|.\label{B8}
\end{eqnarray}

Since $\sum_{i=0}^{(n+1)/2-1}|Y_i|+\sum_{i=(n+1)/2}^n|Y_i|=\sum_{i=0}^n|Y_i|=m$, without loss of generality, assume $\sum_{i=0}^{(n+1)/2-1}|Y_i|\geq\sum_{i=(n+1)/2}^n|Y_i|$.
Let $\sum_{i=(n+1)/2}^n|Y_i|=c$.
Since $m$ is odd, we have $c\leq (m-1)/2$. Then \eqref{B7} is equivalent to
\begin{equation}\label{Z6}
\begin{split}
\textrm{cr}(D^1_0)=&\textrm{cr}(D)+\sum_{i=0}^{(n+1)/2-2}((n+1)/2-i-1)|Y_i|\\&+c+\sum_{i=(n+1)/2}^{n}(i-(n+1)/2)|Y_i|.
\end{split}
\end{equation}
We can obtain a good drawing $D'$ by adjusting the drawing of edge $y_0z_0$ in $N(o,\varepsilon)$ in $D^1_{\frac{n+1}{2}}$ such that
$y_0z_0$ crosses $oy_j$ for $y_j\in \cup_{i=\frac{n+1}{2}}^{n} Y_i$ instead of $oy_j$ for $y_j\in \cup_{i=0}^{\frac{n-1}{2}} Y_i$.
Hence,
\begin{equation}\label{Z7}
\begin{split}
\textrm{cr}(D')=\textrm{cr}(D)+\sum_{i=(n+1)/2}^{n-1}(n-i)|Y_i|+c+\sum_{i=0}^{(n+1)/2-1}i|Y_i|.
\end{split}
\end{equation}
Draw an edge $xo$ next to edge $xy_0$ in $D'$ such that $xo$ crosses an edge $e$ if $xy_0$ crosses $e$ outside $N(o,\varepsilon)$,
which completes the drawing of $D_1$.
For example, a good drawing $D$ of $K_{1,1,5,5}$ with $\sum_{i=0}^{(n+1)/2-1}|Y_i|\geq\sum_{i=(n+1)/2}^n|Y_i|$ is illustrated in Fig.\ref{F5}.
A good drawing $D_1$ constructed from $D$ in Fig.\ref{F5} is illustrated in Fig.\ref{F6}.
Since the number of crossings on $xy_0$ in $D_1$ is $\textrm{cr}_D(E(O,X), E(K_{1,1,m,n})-E(O,X))$, we have
\begin{equation}\label{B9}
\begin{split}
\textrm{cr}(D_1)=\textrm{cr}(D')+\textrm{cr}_D(E(O,X), E(K_{1,1,m,n})-E(O,X)).
\end{split}
\end{equation}

Draw an edge $ox$ next to path $oy_0x$ in $D^1_0$ such that $ox$ crosses $y_0z_i$ in $N(o,\varepsilon)$
for $1\leq i\leq (n-1)/2$ and crosses
an edge if $y_0x$ does outside $N(o,\varepsilon)$,
which completes the drawing of $D_2$.
For example, a good drawing $D_2$ constructed from $D$ in Fig.\ref{F5} is illustrated in Fig.\ref{F7}.
Since the number of crossings on $xy_0$ in $D_2$ is $(n-1)/2+\textrm{cr}_D(E(O,X), E(K_{1,1,m,n})-E(O,X))$, we assert
\begin{equation}\label{B10}
\begin{split}
\textrm{cr}(D_2)=\textrm{cr}(D^1_0)+\textrm{cr}_D(E(O,X), E(K_{1,1,m,n})-E(O,X))+(n-1)/2.
\end{split}
\end{equation}
It is easy to check that both $D_1$ and $D_2$ are good drawings of $K_{1,m+1,n+1}$ with the independent sets $X, \{y_0\}\cup Y$ and $O\cup Z$.
Hence, $\textrm{cr}(D_i)\geq \textrm{cr}(K_{1,m+1,n+1})$ for $1\leq i\leq 2$.
\captionsetup[figure]{labelfont={bf},name={Fig.},labelsep=period}
\begin{figure}[htp]
	\centering
	\begin{minipage}[b]{1.0\textwidth}
		\centering
		\includegraphics[width= 0.85\textwidth]{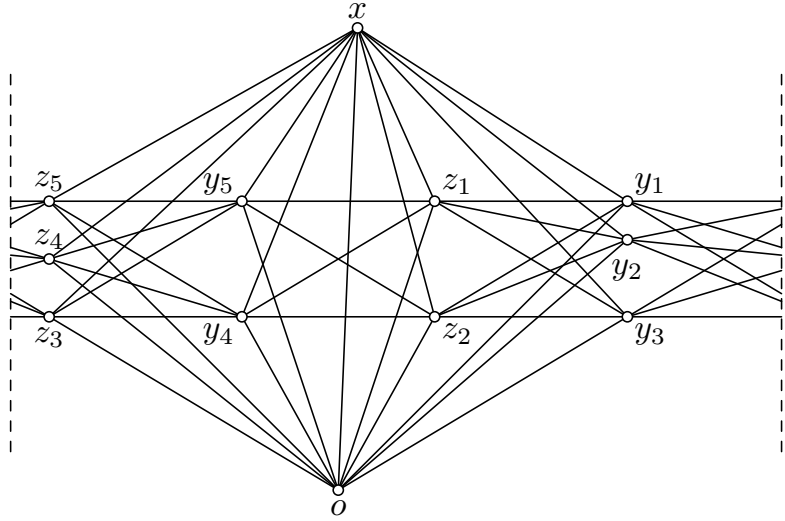}
\begin{minipage}[t]{\textwidth} \caption{A good drawing $D$ of $K_{1,1,m,n}$ in a cylinder with $m=n=5$ and $\sum_{i=0}^{(n+1)/2-1}|Y_i|\geq \sum_{i=(n+1)/2}^n|Y_i|$.}\label{F5} \end{minipage}
	\end{minipage}
\end{figure}
\captionsetup[figure]{labelfont={bf},name={Fig.},labelsep=period}
\begin{figure}[htp]
	\centering
	\begin{minipage}[b]{1.0\textwidth}
		\centering
		\includegraphics[width= 0.85\textwidth]{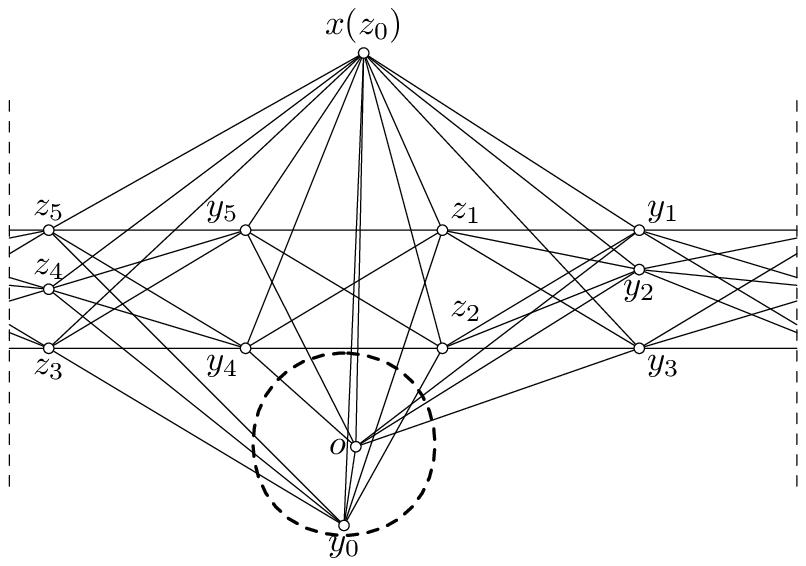}
\begin{minipage}[t]{\textwidth} \caption{A good drawing $D_1$ of $K_{1,6,6}$ in a cylinder constructed from $D$ in Fig.\ref{F5}.
  Both the locations of vertices $o$ and $y_0$ were adjusted for visual clarity.
  The dashed circle denotes the boundary of $N(o,\varepsilon)$ in $D$, which is not a part of $D_1$.}\label{F6} \end{minipage}
	\end{minipage}
\end{figure}
By \eqref{Z6}-\eqref{B10}, we claim
\begin{equation}\label{B11}
\begin{split}
\textrm{cr}(D_1)+\textrm{cr}(D_2)=&2(\textrm{cr}(D)+\textrm{cr}_{D}(E(O,X), E(K_{1,1,m,n})-E(O,X))+c)\\&+\frac{n-1}{2} \sum_{i=0}^n|Y_i|+\frac{n-1}{2}.
\end{split}
\end{equation}
By Lemma \ref{lem3}, $\textrm{cr}_{D}(E(O,X),E(K_{1,1,m,n})-E(O,X))\leq \textrm{cr}(D)-\textrm{cr}(K_{2,m,n}).$
Combining $m=\sum_{i=0}^n|Y_i|$, $c\leq (m-1)/2$ and
$\textrm{cr}(D_i)\geq \textrm{cr}(K_{1,m+1,n+1})$ for $1\leq i\leq 2$,
\eqref{B11} implies
$\textrm{cr}(D)\geq \frac{1}{2}(\textrm{cr}(K_{1,m+1,n+1})+\textrm{cr}(K_{2,m,n})-\frac{1}{4}(m+1)(n+1)+1)$.
By the definition of crossing number, we have $\textrm{cr}(K_{1,1,m,n})\geq\frac{1}{2}(\textrm{cr}(K_{1,m+1,n+1})+\textrm{cr}(K_{2,m,n})-\frac{1}{4}(m+1)(n+1)+1)$.
\end{proof}
\captionsetup[figure]{labelfont={bf},name={Fig.},labelsep=period}
\begin{figure}[htp]
	\centering
	\begin{minipage}[b]{1.0\textwidth}
		\centering
		\includegraphics[width= 1.0\textwidth]{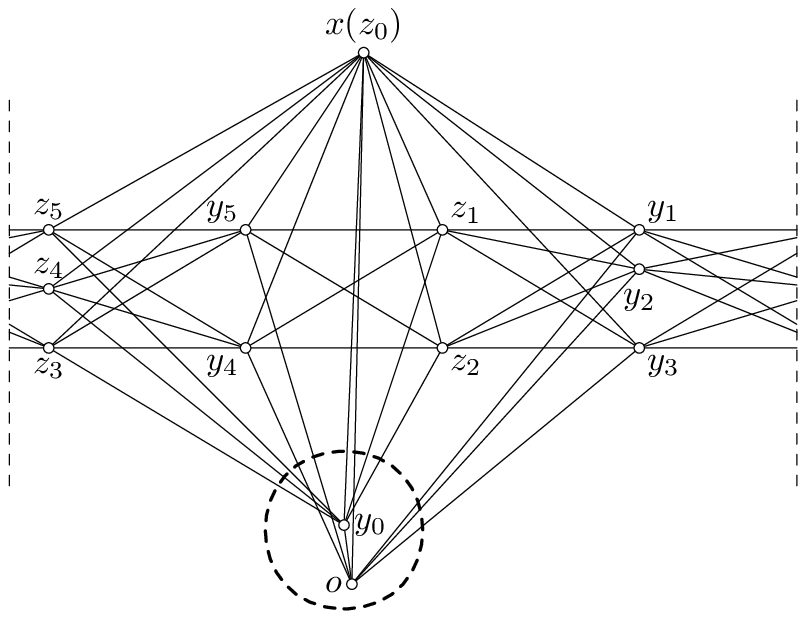}
\begin{minipage}[t]{\textwidth} \caption{A good drawing $D_2$ of $K_{1,6,6}$ in a cylinder constructed from $D$ in Fig.\ref{F5}.
  Both the locations of vertices $o$ and $y_0$ were adjusted for visual clarity.
  The dashed circle denotes the boundary of $N(o,\varepsilon)$ in $D$, which is not a part of $D_2$.}\label{F7} \end{minipage}
	\end{minipage}
\end{figure}
\begin{thm}\label{thm3}
If $m$ is even and $n$ is odd, then
\begin{equation}\nonumber
\begin{split}
\textrm{cr}(K_{1,1,m,n})\geq& \frac{1}{4}(\textrm{cr}(K_{m+1,n+2})+\textrm{cr}(K_{m+3,n+2})+2\textrm{cr}(K_{2,m,n})
\\&-m(n+1)-\frac{1}{4}(n+1)^2).
\end{split}
\end{equation}
\end{thm}

\begin{proof}If $n$ is odd, then repeat the proof we are done in the first paragraph in the proof of Theorem \ref{thm2}.
Let $G_2$ be the graph obtained by adding $ox$ in $G_1$.
Then $G_2$ is isomorphic to $K_{1,m+1,n+1}$ with the independent sets $X, Y\cup \{y_0\}$ and $O\cup Z$.
By drawing an edge $xo$ next to edge $xy_0$ in $D^1_\frac{n+1}{2}$ such that $xo$ crosses an edge outside $N(o,\varepsilon)$
if and only if $xy_0$ does,
we can get a good drawing $D_1$ of $K_{1,m+1,n+1}$.
Hence, $\textrm{cr}(D_1)\geq \textrm{cr}(K_{1,m+1,n+1})$.
For example, a good drawing of $K_{1,1,4,5}$ is illustrated in Fig.\ref{F8},
and a good drawing $D_1$ of $K_{1,5,6}$ constructed from $D$ in Fig.\ref{F8}
is illustrated in Fig.\ref{F9}.
Since the number of crossings on $xy_0$ in $D_1$ is $\textrm{cr}_D(E(O,X), E(K_{1,1,m,n})-E(O,X))$, we have
\begin{equation}\label{J1}
\begin{split}
\textrm{cr}(D_1)=\textrm{cr}(D^1_\frac{n+1}{2})+\textrm{cr}_D(E(O,X), E(K_{1,1,m,n})-E(O,X)).
\end{split}
\end{equation}
\captionsetup[figure]{labelfont={bf},name={Fig.},labelsep=period}
\begin{figure}[htp]
	\centering
	\begin{minipage}[b]{0.9\textwidth}
		\centering
		\includegraphics[width= 1.0\textwidth]{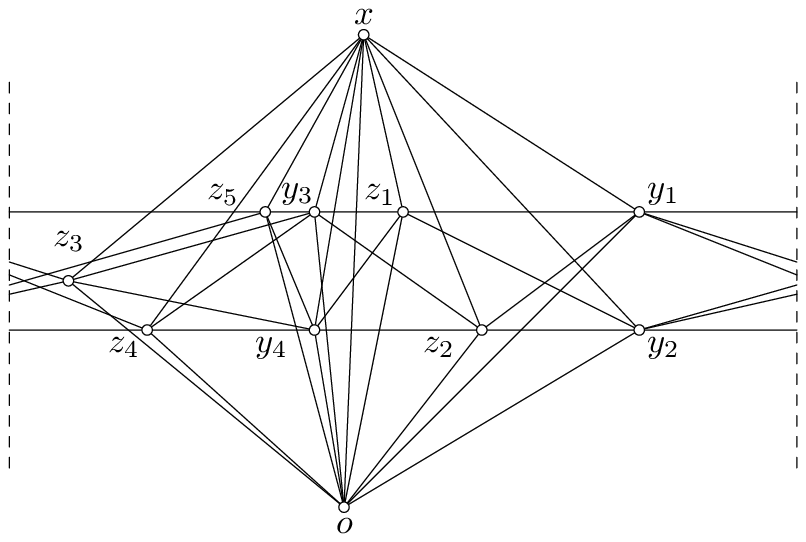}
\begin{minipage}[t]{\textwidth}
\caption{A good drawing $D$ of $K_{1,1,4,5}$ in a cylinder.}\label{F8}
\end{minipage}
	\end{minipage}
	\hfill
	\begin{minipage}[b]{0.9\textwidth}
		\centering
		\includegraphics[width= 1.0\textwidth]{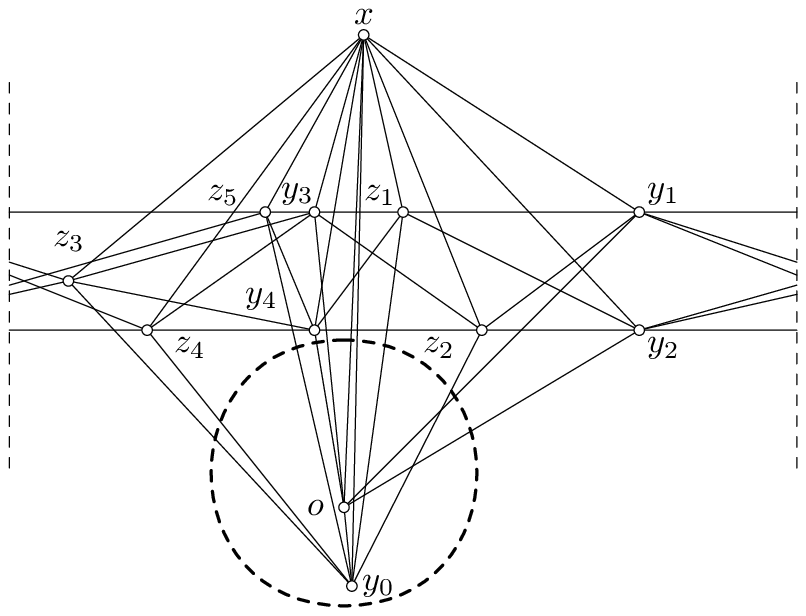}
		\begin{minipage}[t]{\textwidth} \caption{A good drawing $D_1$ of $K_{1,5,6}$ constructed from $D$ in Fig.\ref{F8}.
  Several locations of vertices were adjusted for visual clarity.
  The dashed circle denotes the boundary of $N(o,\varepsilon)$ in $D$, which is not a part of $D_1$.}\label{F9}\end{minipage}
	\end{minipage}
\end{figure}
By deleting $E(X,Z)$ from $D^1_0$, we can get a good drawing $D_2$ of $K_{m+1,n+2}$ with the independent sets $Y\cup\{y_0\}$ and $O\cup X\cup Z$.
For example, a good drawing $D_2$ of $K_{5,7}$ constructed from $D$ in Fig.\ref{F8} is illustrated in Fig.\ref{F10}.
Hence, $\textrm{cr}(D_2)\geq \textrm{cr}(K_{m+1,n+2})$.
Since $\textrm{cr}_{D^1_0}(E(X,Z), E(K_{1,1,m,n})-E(X,Z))=\textrm{cr}_D(E(X,Z),$
$E(K_{1,1,m,n})-E(X,Z))$,
we assert
\begin{equation}\label{J2}
\begin{split}
\textrm{cr}(D_2)=&\textrm{cr}(D^1_0)-\textrm{cr}_{D}(E(X,Z), E(K_{1,1,m,n})-E(X,Z)).
\end{split}
\end{equation}

\captionsetup[figure]{labelfont={bf},name={Fig.},labelsep=period}
\begin{figure}[htp]
		\centering
		\includegraphics[width= 0.85\textwidth]{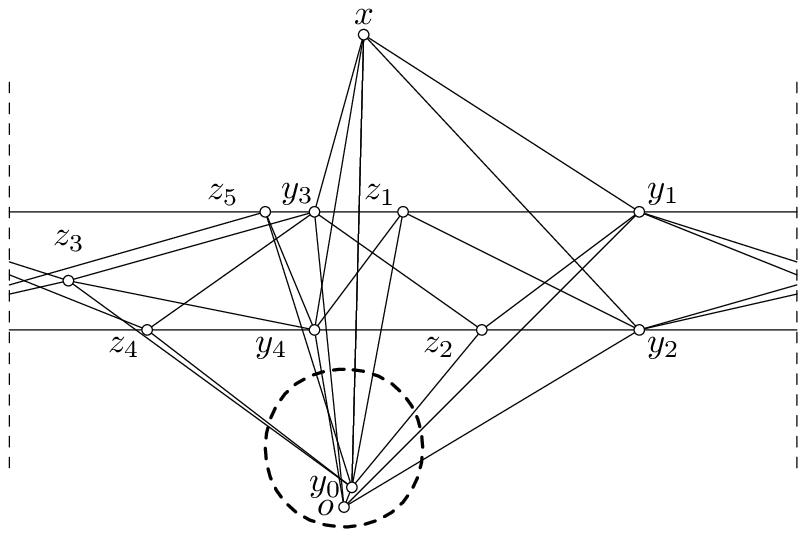}
\caption{A good drawing $D_2$ of $K_{5,7}$ constructed from $D$ in Fig.\ref{F8}.
  Several locations of vertices were adjusted for visual clarity.
  The dashed circle denotes the boundary of $N(o,\varepsilon)$ in $D$, which is not a part of $D_2$.}\label{F10} 
\end{figure}
\captionsetup[figure]{labelfont={bf},name={Fig.},labelsep=period}
\begin{figure}[htp]
		\centering
		\includegraphics[width= 0.85\textwidth]{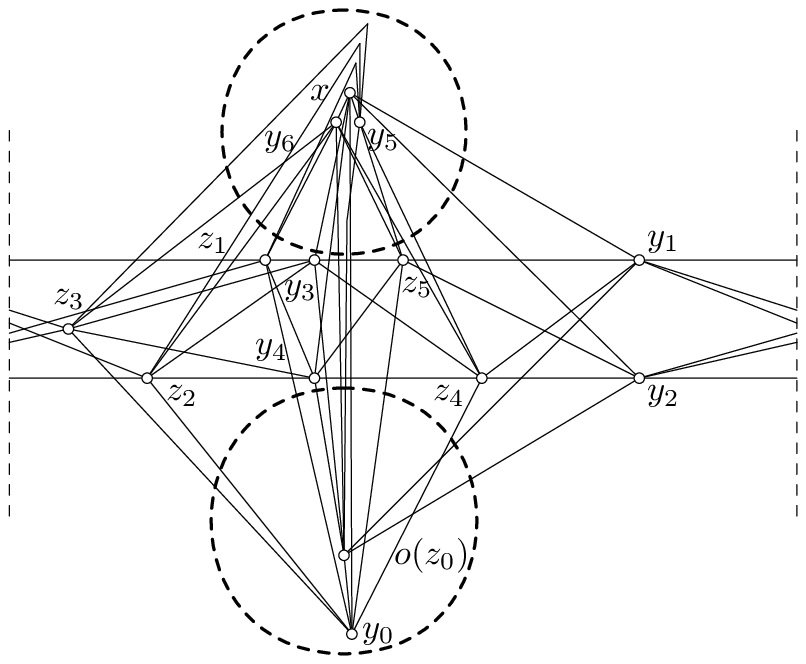}
\caption{A good drawing $D_3$ of $K_{7,7}$ constructed from $D_2$ in Fig.\ref{F9}.
  Several locations of vertices were adjusted for visual clarity.
  The dashed circles denote the boundaries of $N(x,\varepsilon')$ and $N(o,\varepsilon)$ in $D$ respectively, which are not parts of $D_3$.}\label{F11}
\end{figure}
Assume $G_4=G_3^{y_{m+2}Ty_{m+1}}$,
where $V(G_3)=V(G_2)\cup \{y_{m+1}\}$ and $E(G_3)=E(G_2)\cup \{y_{m+1}z_i|1\leq i\leq n\}\cup \{y_{m+1}o,y_{m+1}x\}-\{xz_i|1\leq i\leq n\}\cup \{xo\}$.
Then $G_4$ is isomorphic to $K_{m+3,n+2}$ with the independent sets $\{y_0,y_{m+1},y_{m+2}\}\cup Y$ and $O\cup X\cup Z$.
In $D_1$, we give $o$ another label $z_0$.
By renaming the vertices of $z_i$ if necessary,
we assume $\pi_{D_1}^{O\cup Z}(x)=(z_0,...,z_n)$ and
$Y'_i=\pi_{D_1}(xz_i,xz_{i+1})$ for $0\leq i\leq n$.
Then $m+1=\sum_{i=0}^n|Y'_i|$.
Since $n+1$ is even, we can construct a good drawing $D_3$ of $G_4$ with $\textrm{cr}(D_1)+\textrm{cr}_{D_1}(E(X,O\cup Z),E(G_2)-E(x))+(n+1)/2(m+(n+1)/2)$ crossings by Lemma \ref{lem1}.
Hence, $\textrm{cr}(D_3)\geq \textrm{cr}(K_{m+3,n+2})$.
For example, a good drawing $D_3$ of $K_{7,7}$ constructed from $D_1$ in Fig.\ref{F9} is illustrated in Fig.\ref{F11}.

Since $E(G_2)-E(x)=E(Y\cup \{y_0\},O\cup Z)=E(Y,O\cup Z)\cup E(\{y_0\},O\cup Z)$, by \eqref{I2}, we claim
\begin{equation}\label{X0}
\begin{split}
\textrm{cr}_{D_1}(E(X,O\cup Z),E(G_2)-E(x))=&\textrm{cr}_{D_1}(E(X,O), E(Y,O\cup Z))
\\&+\textrm{cr}_{D_1}(E(X,O), E(\{y_0\},O\cup Z))
\\&+\textrm{cr}_{D_1}(E(X,Z), E(Y,O\cup Z))
\\&+\textrm{cr}_{D_1}(E(X,Z), E(\{y_0\},O\cup Z)).
\end{split}
\end{equation}
By the definition a of good drawing, $\textrm{cr}_{D_1}(E(X,O),E(Y,O))=0$.
By \eqref{I2} and the construction of $D_1$, we have
\begin{equation}\label{X1}
\begin{split}
\textrm{cr}_{D_1}(E(X,O), E(Y,O\cup Z))&=\textrm{cr}_{D_1}(E(X,O),E(Y,Z))\\&=\textrm{cr}_{D}(E(X,O),E(Y,Z)).
\end{split}
\end{equation}
By the construction of $D_1$, we assert
\begin{eqnarray}
\textrm{cr}_{D_1}(E(X,O),E(\{y_0\},O\cup Z))&=&0,\label{X2}\\
\textrm{cr}_{D_1}(E(X,Z), E(Y,O\cup Z))&=&\textrm{cr}_{D}(E(X,Z), E(Y,O\cup Z)),\label{X3}\\
\textrm{cr}_{D_1}(E(X,Z), E(\{y_0\},O\cup Z))&=&\textrm{cr}_{D}(E(X,Z), E(O,Z)).\label{X4}
\end{eqnarray}
By \eqref{I1}, \eqref{I2}, Lemma \ref{lem2} and the definition of a good drawing,
\begin{equation}\label{X5}
\begin{split}
\textrm{cr}_{D}(E(X,Z), E(K_{1,1,m,n})-E(X,Z))=&\textrm{cr}_{D}(E(X,Z), E(O,Z))\\&+\textrm{cr}_{D}(E(X,Z), E(Y,O\cup Z)).
\end{split}
\end{equation}
By \eqref{X0}-\eqref{X5} and Lemma \ref{lem2}, we claim
\begin{equation}\label{C6}
\begin{split}
\textrm{cr}(D_3)=&\textrm{cr}(D_1)+\textrm{cr}_{D}(E(O,X), E(Y,Z))\\&+\textrm{cr}_{D}(E(X,Z), E(K_{1,1,m,n})-E(X,Z))+\frac{n+1}{2}(m+\frac{n+1}{2}).
\end{split}
\end{equation}

By Lemma \ref{lem3}, we have $\textrm{cr}_{D}(E(O,X), E(Y,Z))=\textrm{cr}_{D}(E(O,X), E(K_{1,1,m,n})-E(O,X))$.
By \eqref{J1},\eqref{J2} and \eqref{C6}, we assert
\begin{equation}\label{C7}
\begin{split}
\textrm{cr}(D_2)+\textrm{cr}(D_3)=&2\textrm{cr}_{D}(E(O,X), E(K_{1,1,m,n})-E(O,X))+\textrm{cr}(D^1_0)\\&+\textrm{cr}(D^1_\frac{n+1}{2})+\frac{n+1}{2}(m+\frac{n+1}{2}).
\end{split}
\end{equation}
By putting \eqref{B7} and \eqref{B8} into \eqref{C7}, we claim
\begin{equation}\label{C8}
\begin{split}
\textrm{cr}(D_2)+\textrm{cr}(D_3)=&2(\textrm{cr}(D)+\textrm{cr}_{D}(E(O,X), E(K_{1,1,m,n})-E(O,X)))\\&+\frac{n+1}{2} \sum_{i=0}^{n}|Y_i|+\frac{n+1}{2}(m+\frac{n+1}{2}).
\end{split}
\end{equation}
By Lemma \ref{lem3}, $\textrm{cr}_{D}(E(O,X),E(K_{1,1,m,n})-E(O,X))\leq \textrm{cr}(D)-\textrm{cr}(K_{2,m,n}).$
$m=\sum_{i=0}^n|Y_i|$ and
$\textrm{cr}(D_2)\geq \textrm{cr}(K_{m+1,n+2})$ and $\textrm{cr}(D_3)\geq \textrm{cr}(K_{m+3,n+2})$.
\eqref{C8} implies that
$\textrm{cr}(D)\geq \frac{1}{4}(\textrm{cr}(K_{m+1,n+2})+\textrm{cr}(K_{m+3,n+2})+2\textrm{cr}(K_{2,m,n})-m(n+1)-\frac{1}{4}(n+1)^2)$.
By the definition of crossing number, we have $\textrm{cr}(K_{1,1,m,n})\geq \frac{1}{4}(\textrm{cr}(K_{m+1,n+2})+\textrm{cr}(K_{m+3,n+2})+2\textrm{cr}(K_{2,m,n})-m(n+1)-\frac{1}{4}(n+1)^2)$.
\end{proof}

\begin{cor}\label{cor3}
If at least one of $m$ and $n$ is odd and both ZC and HC on $K_{2,m,n}$ are true, then HC on $K_{1,1,m,n}$ is true.
\end{cor}

\begin{proof} By \eqref{A1}, it suffices to prove $\textrm{cr}(K_{1,1,m,n})\geq Z(m+2,n+2)-mn+\lfloor\frac{m}{2}\rfloor\lfloor\frac{n}{2}\rfloor.$
Since ZC is true, it is easy to check that \eqref{A2} holds by \eqref{A10}, \eqref{A4} and \eqref{A6}.
Since HC on $K_{2,m,n}$ is true, \eqref{A3} holds.

If both $m$ and $n$ are odd, then let $m=2r+1$ and $n=2s+1$.
By Theorem \ref{thm2}, $\textrm{cr}(K_{1,1,2r+1,2s+1})\geq \frac{1}{2}(\textrm{cr}(K_{1,2r+2,2s+2})+\textrm{cr}(K_{2,2r+1,2s+1})-\frac{1}{4}(2r+2)(2s+2)+1)$.
By \eqref{A2}, $\textrm{cr}(K_{1,2r+2,2s+2})=Z(2r+3,2s+3)-(r+1)(s+1)$.
By \eqref{A3}, $\textrm{cr}(K_{2,2r+1,2s+1})=Z(2r+3,2s+3)-(2r+1)(2s+1)$.
Hence, $\textrm{cr}(K_{1,1,2r+1,2s+1})\geq Z(2r+3,2s+3)-(2r+1)(2s+1)+rs=Z(m+2,n+2)-mn+\lfloor\frac{m}{2}\rfloor\lfloor\frac{n}{2}\rfloor$.

Suppose one of $m$ and $n$ is odd and the other is even. Without loss of generality, assume $m$ is even and $n$ is odd.
Let $m=2r$ and $n=2s+1$.
By Theorem \ref{thm3}, $\textrm{cr}(K_{1,1,2r,2s+1})\geq \frac{1}{4}(\textrm{cr}(K_{2r+1,2s+3})+\textrm{cr}(K_{2r+3,2s+3})+2\textrm{cr}(K_{2,2r,2s+1})
-4r(s+1)-(s+1)^2)$.
Since ZC is true, $\textrm{cr}(K_{2r+1,2s+3})=r^2(s+1)^2$ and $\textrm{cr}(K_{2r+3,2s+3})=(r+1)^2(s+1)^2$.
By \eqref{A3}, $\textrm{cr}(K_{2,2r,2s+1})=r(r+1)(s+1)^2-4r(2s+1)$.
Hence, $\textrm{cr}(K_{1,1,2r,2s+1})\geq r(r+1)(s+1)^2-r(3s+2)=Z(m+2,n+2)-mn+\lfloor\frac{m}{2}\rfloor\lfloor\frac{n}{2}\rfloor$.
\end{proof}

\begin{lem}[\cite{Asa}]\label{lem5}
$\textrm{cr}(K_{2,3,n})=Z(5,n+2)-3n$.
\end{lem}

\begin{lem}[\cite{Hua}]\label{lem6}
$\textrm{cr}(K_{1,4,n})=Z(5,n+1)-2\lfloor\frac{n}{2}\rfloor.$
\end{lem}

\begin{cor}\label{cor4}
$\textrm{cr}(K_{1,1,3,n})=Z(5,n)+\lfloor\frac{3n}{2}\rfloor.$
\end{cor}

\begin{proof}By Lemma \ref{lem4}, \ref{lem5}, \ref{lem6}, and the same arguments in the proof of Corollary \ref{cor3},
we have $\textrm{cr}(K_{1,1,3,n})=Z(5,n)+\lfloor\frac{3n}{2}\rfloor$.
\end{proof}

By Corollary \ref{cor1} and \ref{cor3}, we have the following result:

\begin{cor}\label{cor5}
If both ZC and HC on $K_{2,m,n}$ are true, then HC on $K_{1,1,m,n}$ holds.
\end{cor}
\section{Acknowledgements}
The work was supported by the National Natural Science
Foundation of China (No.\,61401186, No.\,11901268).

\end{document}